\documentclass[onecolumn,authoryear]{els-mrw} 

\usepackage{amsmath,amssymb,amsfonts,amsthm,makeidx,graphicx}
\usepackage{txfonts}
\usepackage{helvet}

\newtheorem{ass}[theorem]{Assumption}  

\makeatletter
\DeclareOldFontCommand{\rm}{\normalfont\rmfamily}{\mathrm}
\DeclareOldFontCommand{\sf}{\normalfont\sffamily}{\mathsf}
\DeclareOldFontCommand{\tt}{\normalfont\ttfamily}{\mathtt}
\DeclareOldFontCommand{\bf}{\normalfont\bfseries}{\mathbf}
\DeclareOldFontCommand{\it}{\normalfont\itshape}{\mathit}
\DeclareOldFontCommand{\sl}{\normalfont\slshape}{\@nomath\sl}
\DeclareOldFontCommand{\sc}{\normalfont\scshape}{\@nomath\sc}
\makeatother

\usepackage{csquotes}  
\newcommand\q{\enquote}

\usepackage{physics}   

\newcommand{\tm}{\times}

\DeclareMathOperator*{\esssup}{ess\,sup}

\newcommand \Limsup {\mathop{\overline{\lim}}}

\newcommand \eps {\varepsilon}

\newcommand \N   {\mathbb{N}}
\newcommand \R   {\mathbb{R}}

\newcommand \K   {\mathcal{K}}
\newcommand \Kinf{\mathcal{K_\infty}}

\newcommand \KL  {\mathcal{KL}}

\newcommand \LL  {\mathcal{L}}

\newcommand{\Uc}{\ensuremath{\mathcal{U}}}

\newcommand{\vertiii}[1]{{\left\vert\kern-0.25ex\left\vert\kern-0.25ex\left\vert #1 
    \right\vert\kern-0.25ex\right\vert\kern-0.25ex\right\vert}}

\newcommand \qrq   {\quad\Rightarrow\quad}

\newcommand \qiq   {\quad\Iff\quad}

\newcommand \Iff   {\Leftrightarrow}

\newcommand \id  {\operatorname{id}}

\newcommand{\normt}[1]{{\left\vert\kern-0.25ex\left\vert\kern-0.25ex\left\vert #1 
		\right\vert\kern-0.25ex\right\vert\kern-0.25ex\right\vert}}

\newcommand{\rmd}{\mathrm{d}}

\newif\ifMath					
\newif\ifEngi					

\newif\ifDFGtext					 

\newif\ifAndo              
													
\newif\ifExercises					
\newif\ifSolutions          
\newif\ifGerman							
\newif\ifEnglish						

\newif\ifnothabil						

\newif\ifFuture							

\newif\ifConf                    
\newif\ifJournal								 

\newif\ifNOTFORBOOK
\newif\ifFullVersion
\newif\ifExludedDueToSpaceReasons

\usepackage{xifthen}

\newcommand{\einsnorm}[2]{\ensuremath{
    \!\!\;\!\!\!\;
    \left\bracevert\!\!\!\!\!\left\bracevert
    \!
		\ifthenelse{\isempty{#2}}{#1}{#1(#2)}
    \!
      \right\bracevert\!\!\!\!\!\right\bracevert
    \!\!\;\!\!\!\;
  }}

\usepackage{xcolor}
\definecolor{blond}{rgb}{0.98, 0.94, 0.75}
	
\newlength\mytemplen
\newsavebox\mytempbox

\makeatletter
\newcommand\mybluebox{%
    \@ifnextchar[
       {\@mybluebox}%
       {\@mybluebox[0pt]}}

\def\@mybluebox[#1]{%
    \@ifnextchar[
       {\@@mybluebox[#1]}%
       {\@@mybluebox[#1][0pt]}}

\def\@@mybluebox[#1][#2]#3{
    \sbox\mytempbox{#3}%
    \mytemplen\ht\mytempbox
    \advance\mytemplen #1\relax
    \ht\mytempbox\mytemplen
    \mytemplen\dp\mytempbox
    \advance\mytemplen #2\relax
    \dp\mytempbox\mytemplen
    \colorbox{blond}{\hspace{1em}\usebox{\mytempbox}\hspace{1em}}}

\makeatother

\makeatletter
\let\origd=\d
\renewcommand*\d{
  \relax\ifmmode
    \mathrm{d}%
  \else
    \expandafter\origd
  \fi
}\makeatother

\usepackage{mathtools}

\makeatletter 
\newcommand{\pushright}[1]{\ifmeasuring@#1\else\omit\hfill$\displaystyle#1$\fi\ignorespaces}
\newcommand{\pushleft}[1]{\ifmeasuring@#1\else\omit$\displaystyle#1$\hfill\fi\ignorespaces}
\makeatother 

\newcounter{syscounter}

\newcounter{WPcounter}
\newcounter{PRcounter}

\usepackage{pgfplots}
\usepackage{pgf}
\usepackage{tikz} 
\usetikzlibrary{decorations.pathmorphing} 
\usepackage{tikz-3dplot}
\usepgfplotslibrary{fillbetween}
\usetikzlibrary{arrows}
\usetikzlibrary{graphs,decorations.pathmorphing,decorations.markings}
\usetikzlibrary{calc,math}
\usetikzlibrary{patterns}
\usetikzlibrary{shapes}
\usetikzlibrary{tikzmark}

\usetikzlibrary {positioning}
\usetikzlibrary{shadows}
\usetikzlibrary{patterns.meta}
\usetikzlibrary{plotmarks}

\usepackage{subcaption}  

\usepackage[absolute,overlay]{textpos}
\usepackage{scalefnt}   

\tikzdeclarepattern{
  name=hatch,
  parameters={\hatchsize,\hatchangle,\hatchlinewidth},
  bounding box={(-.1pt,-.1pt) and (\hatchsize+.1pt,\hatchsize+.1pt)},
  tile size={(\hatchsize,\hatchsize)},
  tile transformation={rotate=\hatchangle},
  defaults={
    hatch size/.store in=\hatchsize,hatch size=5pt,
    hatch angle/.store in=\hatchangle,hatch angle=0,
    hatch linewidth/.store in=\hatchlinewidth,hatch linewidth=.4pt,
  },
  code={
      \draw[line width=\hatchlinewidth] (0,0) -- (\hatchsize,\hatchsize);
  }
}

\pgfmathsetmacro\weight{1/2}
\pgfmathsetmacro\third{1/3}
\pgfmathsetmacro\twothirds{2/3}

\tikzset{degil/.style={
            decoration={markings,
            mark= at position 0.5 with {
                  \node[transform shape] (tempnode) {$/$};
                  }
              },
              postaction={decorate}
}
}

\usetikzlibrary{shadows}
\tikzset{
diagonal fill/.style 2 args={fill=#2, path picture={
\fill[#1, sharp corners] (path picture bounding box.south west) -|
                         (path picture bounding box.north east) -- cycle;}},
reversed diagonal fill/.style 2 args={fill=#2, path picture={
\fill[#1, sharp corners] (path picture bounding box.north west) |- 
                         (path picture bounding box.south east) -- cycle;}}
}

\usepgfplotslibrary{fillbetween}
\usetikzlibrary{intersections,through}

\makeatletter 
\tikzset{use path/.code=\tikz@addmode{\pgfsyssoftpath@setcurrentpath#1}}
\makeatother

\definecolor{manipulator-color}{RGB}{88,44,44}
\definecolor{manipulator-contour}{rgb}{0.0, 0.18, 0.39}  

\tikzset{>=latex} 

\begin{document}

\chapter{Input-to-state stability meets small-gain theory}
\label{chap:ISS-for-PDEs}

\author[1]{Andrii Mironchenko}%


\address[1]{\orgname{University of Klagenfurt}, \orgdiv{Department of Mathematics}, \orgaddress{9020, Klagenfurt, Austria}}

\articletag{Chapter Article tagline: update of previous edition,, reprint..}

\maketitle

%
%
%
\begin{abstract}[Abstract]
Input-to-state stability (ISS) unifies global asymptotic stability with respect to variations of initial conditions with robustness with respect to external disturbances. 
First, we present Lyapunov characterizations for input-to-state stability as well as ISS superpositions theorems showing 
relations of ISS to other robust stability properties. 
Next, we present one of the characteristic applications of the ISS framework - the design of event-based control schemes for the stabilization of nonlinear systems.
In the second half of the paper, we focus on small-gain theorems for stability analysis of finite and infinite networks with input-to-state stable components.
First, we present a classical small-gain theorem in terms of trajectories for the feedback interconnection of 2 nonlinear systems. 
Finally, a recent Lyapunov-based small-gain result for a network with infinitely many ISS components is shown.
\end{abstract}

\begin{glossary}[Keywords]
Nonlinear systems; Asymptotic stability; Robustness; Lyapunov functions
\end{glossary}

\section{Introduction}

Consider a system of ordinary differential equations (ODEs) with external inputs of the following form
\begin{subequations}
\label{xdot=f_xu}
\begin{align}
\dot{x} & =  f(x,u),    \label{xdot=f_xu_IntroChap-1}\\
x&(0)  =  x_0. \label{xdot=f_xu_IntroChap-2}
\end{align}
\end{subequations}
Here $x(t)\in\R^n$, $u(t)\in\R^m$, $f:\R^n \times \R^m \to \R^n$ is continuous on $\R^n \times \R^m$, and $x_0\in \R^n$ is a given initial condition. Furthermore, we assume that an input $u$ belongs to the space  $\Uc:= L^{\infty}(\R_+,\R^m)$ of Lebesgue measurable globally essentially bounded functions $u:\R_+ \to \R^m$ endowed with the essential supremum norm 
\begin{eqnarray}
\|u\|_{\infty}:=  \esssup_{t\geq 0}|u(t)| = \inf_{D \subset \R_+,\ \mu(D)=0} \sup_{t \in \R_+ \backslash D} |u(t)|,
\label{eq:Linf_norm}
\end{eqnarray}
where $|x|:=\sqrt{x_1^2 + \ldots + x_k^2}$ denotes the Euclidean norm of a vector $x$. 

We assume that $f$ is Lipschitz continuous with respect to the first argument on bounded subsets, that is, for all $C>0$ there is $L(C)>0$, such that whenever $|x| \leq C$, $|y| \leq C$,  $|v| \leq C$, we have
\begin{eqnarray}
\hspace{-4mm} |f(y,v)-f(x,v)| \leq L(C) |y-x|.
\label{eq:Lipschitz}
\end{eqnarray}
Our assumptions ensure that for each initial condition $x_0$ and each input $u \in\Uc$, the corresponding maximal (Caratheodory) solution of 
\eqref{xdot=f_xu} exists and is unique. We denote it by $\phi(\cdot,x_0,u)$, and it is defined on a certain maximal interval $[0,t_m(x_0,u))$, where $t_m(x_0,u) \in(0,+\infty]$.  
Furthermore, the \emph{boundedness-implies-continuation property} holds, i.e., if $\sup_{t\in[0,t_m(x_0,u))}|\phi(t,x_0,u)|<\infty$, then $t_m(x_0,u)=+\infty$, see \cite[Proposition 1.20]{Mir23}.
If $t_m(x_0,u)=+\infty$ for all $x_0\in\R^n$ and all $u\in\Uc$, then \eqref{xdot=f_xu} is called \emph{forward complete (FC)}.

%
%


To define the stability properties, we will need the following classes of comparison functions.
\begin{equation*}
\begin{array}{ll}
{\K} &:= \left\{\gamma:\R_+\rightarrow\R_+\left|\ \gamma\mbox{ is continuous, strictly} \right. \right. \\
&\phantom{aaaaaaaaaaaaaaaaaaa}\left. \mbox{ increasing and } \gamma(0)=0 \right\}, \\
{\K_{\infty}}&:=\left\{\gamma\in\K\left|\ \gamma\mbox{ is unbounded}\right.\right\},\\
{\LL}&:=\left\{\gamma:\R_+\rightarrow\R_+\left|\ \gamma\mbox{ is continuous and strictly}\right.\right.\\
&\phantom{aaaaaaaaaaaaaaaa} \text{decreasing with } \lim\limits_{t\rightarrow\infty}\gamma(t)=0\},\\
{\KL} &:= \left\{\beta:\R_+\times\R_+\rightarrow\R_+\left|\ \beta \mbox{ is continuous,}\right.\right.\\
&\phantom{aaaaaa}\left.\beta(\cdot,t)\in{\K},\ \beta(r,\cdot)\in {\LL},\ \forall t\geq 0,\ \forall r >0\right\}. \\
\end{array}
\end{equation*}
We refer to \cite[Appendix A]{Mir23} and \cite{Kel14} for a detailed overview of the main properties of these functions.

Let us define the central concept for this paper, introduced in \cite{Son89}, which formalizes the robust stability for nonlinear systems:
\begin{definition}
\label{Def:ISS}
System \eqref{xdot=f_xu} is called \emph{input-to-state stable (ISS)}, if \eqref{xdot=f_xu} is forward complete and there exist $\beta \in \KL$ and $\gamma \in \K$ 
such that for all $x \in \R^n$, all $u\in \Uc$, and all $t\geq 0$ the following holds
\begin {equation}
\label{iss_sum}
|\phi(t,x,u)| \leq \beta(|x|,t) + \gamma( \|u\|_{\infty}).
\end{equation}
\end{definition}
The map $\beta$ describes the transient behavior of the ISS system, and the map $\gamma$ quantifies its asymptotic deviation from the origin. 

Input-to-state stable systems comprise several essential properties. For example, setting $u\equiv 0$ into \eqref{iss_sum}, we see that undisturbed ISS systems are uniformly globally asymptotically stable, i.e., bounded by $\KL$-function $\beta$, which gives a uniform decay rate for solutions starting in any bounded ball. 
This property (under our assumptions on $f$) is equivalent to \emph{global asymptotic stability} (local stability $+$ global attractivity) of the undisturbed system \eqref{iss_sum}, see \cite[Theorem B.37]{Mir23}. 

On the other hand, taking in \eqref{iss_sum} the limit superior when $t\to\infty$, we see that if the system is ISS, then there is $\gamma\in\Kinf$ such that 
\begin{eqnarray}
\limsup_{t\to\infty}|\phi(t,x,u)| \leq  \gamma( \|u\|_{\infty}),\quad x \in\R^n,\quad u\in\Uc.
\label{eq:AG}
\end{eqnarray}
This property is called \emph{asymptotic gain property}, and $\gamma$ is called an \emph{asymptotic gain} of a system \eqref{iss_sum}. 
In particular, the trajectory of any AG system (and, thus, of any ISS system) is bounded if the applied input has a bounded magnitude.

Note that the global asymptotic stability of an undisturbed system and the asymptotic gain property are qualitative properties. They 
tell us that after some time, the trajectories converge to a specific neighborhood of $0$ (or of a ball of radius $\gamma( \|u\|_{\infty})$ around the origin), but these properties do not tell us how fast is this convergence.
In contrast to that, ISS is a quantitative property, which gives precise bounds for convergence times characterized by the transient $\beta$ and asymptotic gain $\gamma$, which is of virtue for the constructive design of robust controllers and observers \cite{KoA01}. 

ISS has not only unified the Lyapunov and input-output stability theories but also allowed for using the synergy of the key tools from Lyapunov stability and input-output stability theories: Lyapunov functions, \q{attractivity+stability} characterizations and small-gain theorems. This revolutionized our view on stabilization of nonlinear systems, design of robust nonlinear observers, stability of nonlinear interconnected control systems, nonlinear detectability theory, and supervisory adaptive control \cite{Mir23,Son08}, and made ISS the dominating stability paradigm in nonlinear control theory, with such diverse applications as robotics, mechatronics, systems biology, electrical and aerospace engineering, to name a few, see \cite{KKK95}. 

Here, we present a compendium of the key results in this theory.
For a comprehensive treatise of the ISS theory for ODE systems, we refer to \cite{Mir23}. 
An overview of the ISS theory of infinite-dimensional systems can be found in \cite{MiP20} and \cite{KaK19}.
For a recent survey on ISS of delay systems, a reader may consult \cite{CKP23}.

\section{ISS superposition theorems}

It is well-known that uniform global asymptotic stability can be characterized as a combination of global attractivity with local stability \cite[Theorem B.37]{Mir23}. Here, we present a corresponding result for ISS.

To start with, every ISS system is \emph{uniformly locally stable (ULS)}, that is, there are $\sigma,\gamma \in\Kinf$ and $r>0$ such that whenever $|x|\leq r$,\ $\|u\|_{\infty} \leq r$,\ and $t \geq 0$, we have
\begin{equation}
\label{GSAbschaetzung}
  \left| \phi(t,x,u) \right| \leq \sigma(|x|) +\gamma(\|u\|_{\infty}).
\end{equation}
To see this, note that for any $\beta\in\KL$ it holds that $\beta(|x|,t)\leq \beta(|x|,0)=:\sigma(|x|)$ for all $x\in\R^n$, $t\geq 0$, and $\sigma \in\Kinf$.
For undisturbed systems, ULS coincides with the classical Lyapunov stability.

Furthermore, every ISS system has the \emph{limit property (LIM)}, namely: there exists
          $\gamma\in\K$ such that 
\begin{equation}
\label{LIM_Absch}
x \in\R^n, \quad  u \in \Uc \qrq \inf_{t \geq 0} |\phi(t,x,u)| \leq \gamma(\|u\|_{\infty}).
\end{equation}    

The following profound result due to \cite[Theorem 1]{SoW96} called the \emph{ISS superposition theorem} states that for forward complete systems, ISS is equivalent to a combination of local stability and a limit property.

\begin{theorem}
\label{thm:ODE-ISS-superposition-theorem}
A system \eqref{xdot=f_xu} is ISS if and only if \eqref{xdot=f_xu} is FC $\wedge$ LIM $\wedge$ ULS.
\end{theorem}

Let us state two important results, which are helpful for the proof of the above superposition theorem and are significant in their own right. The first one, which is due to \cite[Proposition 5.1]{LSW96}, shows that forward complete systems always have bounded finite-time reachability sets.
\begin{proposition}
\label{thm:boundedness_reachability_sets}
\eqref{xdot=f_xu} is forward complete iff it has \emph{bounded reachability sets}, i.e. if for each $r>0$
\[
\sup_{|x|\leq r,\ \|u\|_\infty\leq r,\ t\in[0,r]}|\phi(t,x,u)| < \infty.
\]
\end{proposition}
BRS property implies existence of uniform bounds on finite-time trajectories of families of solutions starting in bounded balls. 
This makes BRS a bridge between the pure well-posedness theory, which studies existence and uniqueness properties of solutions, and the stability theory, which is interested in global in time bounds on families of solutions. 

For systems without inputs, the equivalence between FC and BRS is clear by compactness of bounded balls in finite-dimensional Euclidean spaces. However, for systems with inputs, the set 
$\{(t,x,u):|x|\leq r, \|u\|_\infty\leq r, t\in[0,r]\}$ is not compact, which makes the proof of Theorem~\ref{thm:ODE-ISS-superposition-theorem} much more involved, see \cite{LSW96}. 

The next property introduced in \cite{MiW18b} is a uniform version of the LIM property.
\begin{definition}
\label{def:Limit-properties}
We say that a forward complete system \eqref{xdot=f_xu} has  \emph{uniform limit (ULIM) property}, if there exists
    $\gamma\in\K\cup\{0\}$ so that for every $\eps>0$ and for every $r>0$ there
    exists a $\tau = \tau(\eps,r)$ such that 
for all $x\in B_r$ and all $u\in B_{r,\Uc}$ there is a $t\leq\tau$ such that 
\begin{eqnarray}
|\phi(t,x,u)| \leq \eps + \gamma(\|u\|_{\infty}).
\label{eq:sLIM_ISS_section}
\end{eqnarray}
\end{definition}

It turns out that for ODE systems, the limit property is intrinsically uniform, as was shown in \cite[Proposition 13]{MiW18b} based on 
\cite[Corollary III.3]{SoW96}.
\begin{proposition}
\label{prop:ULIM_equals_LIM_in_finite_dimensions}
\eqref{xdot=f_xu} is LIM $\qiq$ \eqref{xdot=f_xu} is ULIM.
\end{proposition}

Propositions~\ref{thm:boundedness_reachability_sets}, \ref{prop:ULIM_equals_LIM_in_finite_dimensions} 
exploit significantly the properties of ODE systems and do not hold in general for infinite-dimensional systems; see the counterexamples for time-delay systems \cite{MaH23,CWM24} and for infinite networks of ODE systems \cite[Theorem 5]{MiW18b}. 

With Propositions~\ref{thm:boundedness_reachability_sets}, \ref{prop:ULIM_equals_LIM_in_finite_dimensions} in mind, one can derive 
 Theorem~\ref{thm:ODE-ISS-superposition-theorem} from general ISS superposition theorem for infinite-dimensional systems shown in \cite{MiW18b} stating that ISS is equivalent to a combination of BRS $\wedge$ ULIM $\wedge$ ULS. See \cite[Section 2.5]{Mir23} for the detailed proof.

\section{Lyapunov functions}

Lyapunov functions are an indispensable tool for the verification of ISS for nonlinear systems:
\begin{definition}
\index{ISS Lyapunov function}
\index{ISS!Lyapunov function (dissipative form)}
\label{def:ISS_LF-dissipative}
A continuous function $V:\R^n \to \R_+$ is called an \emph{ISS Lyapunov function in a dissipative form} for \eqref{xdot=f_xu},  if there exist
$\psi_1,\psi_2 \in \Kinf$, $\alpha \in \Kinf$ and $\xi \in \K$ 
such that 
\begin{equation}
\label{LyapFunk_1Eig}
\psi_1(|x|) \leq V(x) \leq \psi_2(|x|) \quad \forall x \in \R^n,
\end{equation}
and for any $x \in \R^n, u \in \Uc$ the following inequality holds:
\begin{eqnarray}
\dot{V}_u(x) \leq -\alpha(V(x)) + \xi(\|u\|_\infty),
\label{eq:ISS-LF-DissipativeIneq}
\end{eqnarray}
where \emph{Lie derivative} $\dot{V}_u(x)$ corresponding to the pair $(x,u)$, is the upper right-hand Dini derivative at zero of the function $t\mapsto V(\phi(t,x,u))$, that is:
\begin{equation}
\label{LyapAbleitung-inputs}
\dot{V}_u(x):=\Limsup \limits_{t \rightarrow +0} {\frac{1}{t}\big(V(\phi(t,x,u))-V(x)\big) }.
\end{equation}
\index{derivative!Lie}
\index{Lie derivative}
\end{definition}
The dissipative form of a Lyapunov function resembles the storage function in a theory of dissipative systems, which explains the name \q{dissipative}. 
The following fundamental result is due to \cite{SoW95}:

\begin{theorem}[Smooth converse ISS Lyapunov theorem]
\label{thm:ISS_Converse_Lyapunov_Theorem}
Let $f$ be Lipschitz continuous on bounded balls of $\R^n \tm \R^m$, and $f(0,0)=0$. 
Then \eqref{xdot=f_xu} is ISS if and only if there is an infinitely differentiable ISS Lyapunov function for \eqref{xdot=f_xu}.
\end{theorem}
To show ISS, the existence of a continuous ISS Lyapunov function suffices, see \cite[Theorem 2.12]{Mir23}. The proof of the converse implication is much more challenging and exploits a converse Lyapunov theorem for robust asymptotic stability of nonlinear ODE systems \cite{LSW96}. 
The smoothing technique used in \cite{LSW96} is distinctly finite-dimensional and cannot be generalized to general infinite-dimensional systems. However, if we aim for a merely Lipschitz continuous Lyapunov function, then a general converse Lyapunov theorem for robust asymptotic stability of nonlinear infinite-dimensional systems \cite[Chapter 3]{KaJ11b} can be used to obtain infinite-dimensional extensions of Theorem~\ref{thm:ISS_Converse_Lyapunov_Theorem}, see \cite{MiW17c}.

\section{Lyapunov-based design of event-based controllers}
\label{sec:Event-based control}

Even though the stabilization of nonlinear control systems is already a challenging problem, in real-world applications of control theory, many additional complications appear that we summarized in Table~\ref{tab:Obstructions-nonlinear-control} (taken from \cite[Table 5.1]{Mir23}).
ISS framework can systematically handle these challenges; see \cite[Section 5.11.2]{Mir23} for an overview.

\begin{table}[ht]
\centering
\begin{tabular}{|l|l|}
\hline
\qquad\qquad\qquad\qquad \textbf{Problem} & \qquad\qquad\qquad \textbf{Remedy} \\ 
\hline
Signals are transmitted at discrete moments of time & Sampled-data control,  \\ 
																															 & Event-based control \\ \hline
The full state cannot be measured                      & Observer/Estimator design      \\ \hline
Limitations of the transmission capacity of  			     & Theory of networked control systems \\ 
the communication network 									 			     &  \\ \hline
Delays during the transmission                         & Delay compensation technique   \\ \hline
Unmodeled dynamics                                     & Robust control, ISS  \\ \hline
Measurement errors, actuator errors                    & Robust observers/controllers, ISS  \\ \hline
Quantization of the sensor and input signals           & Quantized controllers \\ \hline
Boundedness of the physical input                      & Controllers with saturation \\ \hline
Uncertainties of parameters                            & Adaptive control     \\ \hline
\end{tabular}
\caption{Various challenges, appearing in real-world applications of control theory and the approaches how to overcome them.}
\label{tab:Obstructions-nonlinear-control}
\end{table}

Typically, one assumes that \q{idealistic} controller (in the absence of all these obstructions) has already been constructed, and the difference between the ideal trajectory and real trajectory is understood as the disturbance. Then, one exploits the ISS properties w.r.t. these disturbances to redesign the idealistic controller to achieve stabilization in a realistic scenario.
In this section, we present this scheme by discussing the ISS-based design of event-based controllers for nonlinear systems due to \cite{Tab07}.

In classical stabilization theory, one assumes that we can update a controller continuously. 
However, digital controllers do not satisfy this assumption, since the controller can be updated only at discrete moments. Hence, a question arises whether the constructed \q{continuous-time controller} will be able to stabilize the system and how the triggering times should be chosen.

Consider a control system
\begin{eqnarray}
\dot{x} = f(x,u),
\label{eq:Event-based-control}
\end{eqnarray}
with $f$ that is locally Lipschitz continuous in both arguments.
Assume that a Lipschitz continuous on bounded balls feedback controller 
\begin{eqnarray}
u:=k(x),
\label{eq:Feedback-controller-event-based}
\end{eqnarray}
is given, making the closed-loop system 
\begin{eqnarray}
\dot{x} = f(x,k(x+e)),
\label{eq:Event-based-control-closed-loop}
\end{eqnarray}
ISS with respect to the \emph{measurement error} $e$.
\index{measurement error}

Since $f$ is Lipschitz continuous on bounded balls of $\R^n \tm \R^m$ and $k$ is Lipschitz continuous on bounded balls, then $(x,e)\mapsto f(x,k(x+e))$ is Lipschitz continuous on bounded balls in both $x$ and $e$.
This implies the well-posedness of \eqref{eq:Event-based-control-closed-loop}. Moreover, Theorem~\ref{thm:ISS_Converse_Lyapunov_Theorem} ensures that there is a smooth ISS Lyapunov function $V:\R^n \to \R_+$ in a dissipative form for \eqref{eq:Event-based-control-closed-loop}. This means that there are $\psi_1,\psi_2 \in \Kinf$, $\alpha \in \Kinf$ and $\xi \in \K$ such that 
\begin{equation}
\label{LyapFunk_1Eig-event}
\psi_1(|x|) \leq V(x) \leq \psi_2(|x|) \quad \forall x \in \R^n,
\end{equation}
and for all $x \in \R^n, e \in \R^n$ the following inequality holds:
\begin{eqnarray}
\nabla V(x) \cdot f(x,k(x+e)) \leq -\alpha(|x|) + \xi (|e|).
\label{eq:ISS-LF-DissipativeIneq-event-based}
\end{eqnarray}

For the \emph{digital implementation} of the controller \eqref{eq:Feedback-controller-event-based},
 we employ the \emph{sample-and-hold technique}, according to which, one performs the measurements $(x(t_k))_{k\in\N}$ of the state at time instants $(t_k)_{k\in\N}$ and employs the corresponding \q{digital} control law
\begin{eqnarray}
u_d(t):= u(t_k)=k(x(t_k)),\quad t \in [t_k,t_{k+1}).
\label{eq:digitalized-control-event-based}
\end{eqnarray}
We are going to update our controller when the error becomes \q{too large}. 
Let us introduce the \emph{error function}
\begin{eqnarray}
e(t):=x(t_k)-x(t),\quad t\in [t_k,t_{k+1}).
\label{eq:Error-event-based-control}
\end{eqnarray}
The system \eqref{eq:Event-based-control} with the digital controller \eqref{eq:digitalized-control-event-based} becomes
\begin{eqnarray}
\dot{x}(t) = f\big(x(t),k\big(x(t_k)\big)\big),\quad t \in [t_k,t_{k+1}),
\label{eq:NL-sys-with-event-based controller}
\end{eqnarray}
which is \eqref{eq:Event-based-control-closed-loop} with the error $e$ as in \eqref{eq:Error-event-based-control}.

In view of \eqref{eq:ISS-LF-DissipativeIneq-event-based}, if for some $\sigma\in(0,1)$
\begin{eqnarray}
\xi (|e|)\leq \sigma \alpha(|x|),
\label{eq:small-error-condition}
\end{eqnarray}
then 
\begin{eqnarray}
\nabla V(x) \cdot f\big(x,k(x+e)\big) \leq -(1-\sigma)\alpha(|x|).
\label{eq:ISS-LF-DissipativeIneq-event-based-decay}
\end{eqnarray}
To ensure that \q{small error condition} \eqref{eq:small-error-condition} always holds, it is enough to update the controller whenever the following condition (or \emph{event}) holds:
\begin{eqnarray}
\xi (|e|)\geq \sigma \alpha(|x|).
\label{eq:big-error-condition}
\end{eqnarray}
Under the event-triggering condition \eqref{eq:big-error-condition}, the estimate \eqref{eq:ISS-LF-DissipativeIneq-event-based-decay} holds along all trajectories of \eqref{eq:NL-sys-with-event-based controller}, $V$ is a strict Lyapunov function for \eqref{eq:NL-sys-with-event-based controller} and thus \eqref{eq:NL-sys-with-event-based controller} is UGAS.

The above analysis makes sense if the event-triggered control is \emph{feasible}, i.e., the sequence of triggering times $(t_k)$ does not have a \emph{Zeno behavior} (i.e., provided that it is unbounded).
The following theorem ensures this under the assumption that $\alpha^{-1}$ and $\xi$ are locally Lipschitz on bounded sets.
\begin{theorem}
\label{thm:Event-based-controller} 
Consider a system \eqref{eq:Event-based-control} with $f$ being Lipschitz continuous on bounded balls of $\R^n \tm \R^m$, and with a locally Lipschitz controller given by \eqref{eq:Feedback-controller-event-based} making the closed-loop system \eqref{eq:Event-based-control-closed-loop} ISS w.r.t. the measurement disturbance $e$. 

Let $V$ be an ISS Lyapunov function for the closed-loop system as defined above with $\alpha^{-1}$ and $\xi$ being locally Lipschitz on bounded sets.

Then for every compact set $S \subset \R^n$ with $0\in S$ there is a time $\tau>0$ such 
 that for any initial condition in $S$, the inter-execution times 
$(t_{i+1}-t_i)_{i\in\N}$ implicitly defined by the execution rule \eqref{eq:big-error-condition} with $\sigma\in(0,1)$
are lower bounded by $\tau$, i.e., $t_{i+1}-t_i\geq \tau$ for any $i\in\N$.
\end{theorem}

Theorem~\ref{thm:Event-based-controller} is due to \cite{Tab07}, where it was
proved for a more general scenario where the update of the control
law is delayed by some positive time $\delta \ge 0$ representing the time required to read
the state from the sensors, compute the value of the controller, and update the actuators.

For a systematic account of the control applications of the ISS framework, we refer to \cite[Chapter 5]{Mir23}.

\section{Small-gain theorems for 2 systems}

Analyzing and controlling nonlinear systems is a complex endeavor, particularly for large-scale systems. The lack of general methods to construct ISS Lyapunov functions for nonlinear systems means that direct stability analysis via construction of an ISS Lyapunov function, is rarely feasible for large-scale systems. In the following sections, we demonstrate powerful nonlinear small-gain theorems that ensure the input-to-state stability of large-scale networks composed of input-to-state stable nonlinear components, given that the network's interconnection structure satisfies the small-gain condition. The small-gain approach has broad applications, including designing robust controllers and observers for nonlinear systems, developing decentralized observers for large-scale networks, and achieving synchronization in multiagent systems.

Consider a feedback interconnection of two systems of the form
\begin{eqnarray}
\label{ISSGesamtNsysteme}
\dot{x}_{1}&=&f_{1}(x_{1},x_{2},u),\\
\dot{x}_{2}&=&f_{2}(x_{1},x_{2},u).
\end{eqnarray}
Here $u \in \Uc:=L^{\infty}(\R_+,\R^m)$ is an \emph{external input} that we, without loss of generality, assume to be the same for all subsystems. Otherwise, we can collect the inputs for individual subsystems $u_1,u_2$ into a vector $(u_1,u_2)$ that will be a common input for all subsystems. The state of $x_i$-subsystem belongs to the space $X_i:=\R^{N_i}$, for some $N_i \in\N$.
We denote the dimension of the state space of the whole system as $N:=N_1 +N_2$. 

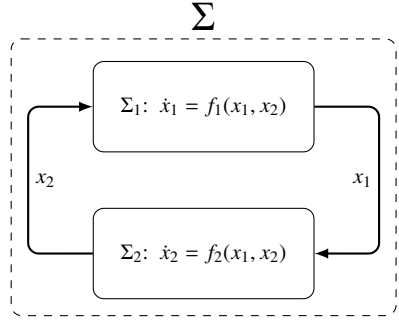
\begin{figure}
\centering
\begin{tikzpicture}[minimum width=3.9cm, text centered, rounded corners, minimum height=1.6cm, scale = 0.75, transform shape]
\node [rectangle, draw](B1)   at (0,0) {\large $\Sigma_1$:\ \ $\dot{x}_1=f_1(x_1,x_2)$};
\node [rectangle, draw](B2)   at (0,-2.6) { \large $\Sigma_2$:\ \ $\dot{x}_2=f_2(x_1,x_2)$};
\draw[->, thick] (B1) to (3.1,0) to (3.1,-2.6) to (B2);
\draw[->, thick] (B2) to (-3.1,-2.6) to (-3.1,0) to (B1);

\node (B3)   at (2.8,-1.3) {\large $x_1$};
\node (B4)   at (-2.8,-1.3){\large $x_2$};
\node (B4)   at (0,1.6) { \huge $\Sigma$};



\draw [rounded corners,dashed] (-3.4,1.2)--(3.4,1.2)--(3.4,-3.7)--(-3.4,-3.7)--cycle;

\end{tikzpicture}
\caption{Feedback interconnection}
\end{figure}

Uniform global asymptotic stability of both $x_1$-subsystem and $x_2$-subsystem in absence of inputs (i.e., if $\dot{x}_1 = f_1(x_1,0,0)$ and $\dot{x}_2 = f_2(0,x_2,0)$ are UGAS), does not guarantee any global properties of the coupled system \eqref{ISSGesamtNsysteme}. 
E.g., consider a coupled system
\begin{align*}
\dot{x_1} &= -x_1 + x_2x_1^2,\\ 
\dot{x_2} &= -x_2. 
\end{align*}
Both subsystems of this system are UGAS in the absence of inputs, and the second subsystem does not have any input. Nevertheless, the coupled system is not forward complete (and thus, it is not ISS), which is easy to see by choosing $x_2$ and $x_1$ large enough. 

This means that for the stability analysis of interconnected systems, we should go beyond the classical dynamical systems paradigm and exploit the information about the influence of subsystems on each other. The following result, which is a special case of \cite[Theorem 2.1]{JTP94}, plays a crucial role in stability analysis of interconnections:

For an $L^\infty$-function $z$ we denote by $\left\|z\right\|_{[0,t]}$ the $L^\infty$-norm of $z$ on the interval $[0,t]$.

\begin{theorem}
\label{thm:SGT-trajectories} 
Assume that both subsystems of a system \eqref{ISSGesamtNsysteme} are ISS. More precisely, we assume that there exist $\beta_1,\beta_2 \in\KL$, internal gains $\gamma_{12},\gamma_{21}$, and external gains $\gamma_1,\gamma_2$ such that for all initial states, and all internal and external $L^\infty$-inputs, as well as for all times $t\geq 0$ the following inequalities hold:
\begin{eqnarray}
\label{eq:ISS_n_sys_sum}
|\phi_1(t,x_2,u)|  &\leq&  \beta_1\left(\left|x_1\right|,t\right) + \gamma_{12}\left(\left\|x_2\right\|_{[0,t]}\right) + \gamma_1\left(\left\|u\right\|_{\infty}\right),\\
|\phi_2(t,x_1,u)|  &\leq&  \beta_1\left(\left|x_2\right|,t\right) + \gamma_{12}\left(\left\|x_1\right\|_{[0,t]}\right) + \gamma_2\left(\left\|u\right\|_{\infty}\right).
\end{eqnarray}
Assume that there is $\rho\in\Kinf$ such that the following small-gain condition holds:
\begin{eqnarray}
(\id+\rho) \circ \gamma_{12} \circ(\id+\rho) \circ\gamma_{21}(r)<r \quad \forall r>0.
\label{eq:Strong-small-gain-condition-cyclic}
\end{eqnarray}
Then, the interconnection \eqref{ISSGesamtNsysteme} is ISS.
\end{theorem}

Define the \emph{gain operator} $\Gamma:\R^2_+\to\R^2_+$:
\begin{eqnarray}
\Gamma(s):=
\begin{pmatrix}
\gamma_{12}(s_2)\\
\gamma_{21}(s_1)
\end{pmatrix}
,\quad s=(s_1,s_2)\in\R^2_+.
\end{eqnarray}
We say that $a \not\geq b$ for some $a=(a_1,a_2),b=(b_1,b_2)\in\R^2_+$  if $a_j < b_j$ for some $j\in\{1,2\}$. 
Now the condition \eqref{eq:Strong-small-gain-condition-cyclic} can be equivalently written in the following form:
\begin{equation}
\label{eq:Strong_SGC-chap-networks}
  (\id + \rho)\circ \Gamma(s) \not\geq s,\quad s\in \R^2_+\setminus\{0\}.
\end{equation}

In \cite{DRW07}, a small-gain theorem for feedback interconnections of an arbitrary finite number of nonlinear ODE systems has been shown by introducing the gain operator for such couplings and exploiting the small-gain condition in the form \eqref{eq:Strong_SGC-chap-networks} 
to guarantee the stability of the whole interconnection.
The proof of Theorem \ref{thm:SGT-trajectories} and its extensions in \cite{DRW07} is based on the ISS superposition theorem (Theorem~\ref{thm:ODE-ISS-superposition-theorem}).

A Lyapunov-based version of the small-gain theorem for feedback coupling of 2 systems has been proposed in \cite{JMW96}, and its extension to couplings of $n$ ODE systems has been proposed in \cite{DRW10}.

Small-gain theorems play an essential role in the analysis of networks with ISS components. They can be used to control coupled systems in combination with the gain assignment technique proposed in \cite[Theorem 2.2]{JTP94}.

\section{Nonlinear small-gain theorem for infinite networks}

We are surrounded by various types of networks, including social networks, power grids, transportation, and manufacturing networks. These networks expand from year to year, and new technologies like cloud computing and fifth-generation (5G) communication further highlight this trend. Since the stability of networks can deteriorate as the number of participating agents grows, it is natural to study the stability of infinite networks, which serve as an overapproximation for large-scale networks. In this section, we discuss a recent generalization of the nonlinear small-gain theorem to interconnections of a countable number of systems.

%

Consider a family of control systems of the form
\begin{equation}\label{eq_subsystem}
  \Sigma_i:\quad \dot{x}_i = f_i(x_i,\bar{x}_i,u_i),\quad i \in \N.%
\end{equation}
This family comes with sequences $(n_i)_{i\in\N}$ and $(m_i)_{i\in\N}$ of positive integers as well as \emph{finite} (possibly empty) sets $I_i \subset \N$, $i \notin I_i$, such that the following holds:
\begin{itemize}
\item The \emph{state vector} $x_i$ is an element of $\R^{n_i}$.%
\item The \emph{internal input vector} $\bar{x}_i$ is composed of the state vectors $x_j$, $j \in I_i$, and thus is an element of $\R^{N_i}$, where $N_i := \sum_{j \in I_i}n_j$.%
\item The \emph{external input vector} $u_i$ is an element of $\R^{m_i}$.%
\item The \emph{right-hand side} $f_i:\R^{n_i} \tm \R^{N_i} \tm \R^{m_i} \rightarrow \R^{n_i}$ is a continuous function.%
\item For every initial state $x_{i0} \in \R^{n_i}$ and all essentially bounded inputs $\bar{x}_i(\cdot)$ and $u_i(\cdot)$, there is a unique solution of $\Sigma_i$, which we denote by $\phi_i(t,x_{i0},\bar{x}_i,u_i)$ (it may be defined only on a bounded time interval).%
\end{itemize}
For each $i\in\N$, we fix norms on the spaces $\R^{n_i}$ and $\R^{m_i}$, respectively (these norms can be chosen arbitrarily). For brevity in notation, we avoid adding an index to these norms, indicating to which space they belong, and write $|\cdot|$ for each of them. The interconnection of the systems $\Sigma_i$, $i\in\N$, is defined on the state space $X := \ell_{\infty}(\N,(n_i))$, where%
\begin{align*}
  \ell_{\infty}(\N,(n_i)) := \{ x = (x_i)_{i\in\N} : x_i \in \R^{n_i},\ \sup_{i \in \N}|x_i| < \infty \},
\end{align*}
which is a Banach space with the $\ell_{\infty}$-type norm%
\begin{equation*}
  \|x\|_X := \sup_{i\in\N}|x_i|.%
\end{equation*}
Similarly, we define the space of admissible external input values:
\begin{equation*}
 U := \ell_{\infty}(\N,(m_i)),\quad \|u\|_U := \sup_{i\in\N}|u_i|.%
\end{equation*}
We choose the class of admissible external input functions as
\begin{equation*}
  \Uc := \{u \in L^{\infty}(\R_+,U) : u \mbox{ is piecewise right-continuous} \},%
\end{equation*}
equipped with the $L^{\infty}$-norm%
\begin{equation*}
  \|u\|_{\Uc} := \esssup_{t\in\R_+} \|u(t)\|_U.%
\end{equation*}
Define the right-hand side of the interconnected system by
\begin{equation*}
  f:X \tm U \rightarrow \prod_{i\in\N}\R^{n_i},\quad  f(x,u) := (f_i(x_i,\bar{x}_i,u_i))_{i\in\N}.%
\end{equation*}
Now, we can represent the interconnected system as follows
\begin{equation*}
  \Sigma:\quad \dot{x} = f(x,u).%
\end{equation*}
To fully define the system, we need to introduce an appropriate notion of solution. For a fixed $(u,x^0) \in \Uc \tm X$, we call a function $\lambda:J \rightarrow X$, where $J = [0,T) \subset \R$ with $0 < T \leq \infty$, a \emph{solution} of the Cauchy problem%
\begin{equation*}
  \dot{x} = f(x,u),\quad x(0) = x^0,%
\end{equation*}
provided that $s \mapsto f(\lambda(s),u(s))$ is a locally integrable $X$-valued function (in the Bochner integral sense) and%
\begin{equation*}
  \lambda(t) = x^0 + \int_0^t f(\lambda(s),u(s))\, \rmd s \mbox{\quad for all\ } t \in J.%
\end{equation*}

\begin{ass}\label{ass:Well-posedness-BIC}
We assume that the system $\Sigma$ is \emph{well-posed}, that is, for every initial state $x^0 \in X$ and every input $u \in \Uc$, there exists a unique maximal solution. We denote this solution by $\phi(\cdot,x^0,u):[0,t_{\max}(x^0,u)) \rightarrow X$, where $0 < t_{\max}(x^0,u) \leq \infty$.

Furthermore, let all uniformly bounded maximal solutions $\phi(\cdot,x,u)$ of $\Sigma$ exist on $\R_+$.
\end{ass}

The construction of an ISS Lyapunov function is a complex problem, which becomes especially challenging if the system is nonlinear and of large size. In this paper, we assume that all components $\Sigma_i$ of an infinite network are ISS with corresponding ISS Lyapunov functions $V_i$. To find an ISS Lyapunov function $V$ for $\Sigma$, we exploit the interconnection structure and construct $V$ from $V_i$. Hence, we make the following assumption.

\begin{ass}\label{ass_subsystem_iss}
For each $i \in \N$, there exists a continuous function $V_i:\R^{n_i} \rightarrow \R_+$ which is continuously differentiable outside of $x_i = 0$ and satisfies the following properties:%
\begin{enumerate}
\item[(L1)] There exist $\psi_{1},\psi_{2} \in \Kinf$ (same for all $i$) such that%
\begin{equation}\label{eq_subsystem_iss_coerc} 
  \psi_{1}(|x_i|) \leq V_i(x_i) \leq \psi_{2}(|x_i|) \mbox{\quad for all\ } x_i \in \R^{n_i}.%
\end{equation}
\item[(L2)] There exist $\gamma_{ij} \in \K \cup \{0\}$, where $\gamma_{ij} = 0$ for all $j \in \N \setminus I_i$, and $\gamma_{iu} \in \K$ as well as $\tilde\alpha \in \K$ (same for all $i$) such that for all $x = (x_j)_{j\in\N} \in X$ and $u = (u_j)_{j\in\N} \in U$ the following implication holds:%
\begin{align}\label{eq_subsystem_orbitalder_est}
\begin{split}
  V_i(x_i) &> \max\Bigl\{\sup_{j\in I_i}\gamma_{ij}(V_j(x_j)),\gamma_{iu}(|u_i|) \Bigr\} \\
	&\Rightarrow \nabla V_i(x_i)f_i(x_i,\bar{x}_i,u_i) \leq -\tilde\alpha(V_i(x_i)).%
\end{split}
\end{align}
\end{enumerate}
The function $V_i$ is called an \emph{ISS Lyapunov function for $\Sigma_i$}. The functions $\gamma_{ij}$ and $\gamma_{iu}$ are called \emph{internal gains} and \emph{external gains}, respectively.%
\end{ass}

Using the internal gains $\gamma_{ij}$ from Assumption \ref{ass_subsystem_iss}, we define the \emph{gain operator} $\Gamma:\ell_{\infty}^+ \rightarrow \ell_{\infty}^+$ by
\begin{equation}
\label{eq_gamma_def}
  \Gamma(s) := \Bigl( \sup_{j\in\N}\gamma_{ij}(s_j) \Bigr)_{i\in\N}.%
\end{equation}
The following assumption guarantees that $\Gamma$ is well-defined and continuous.
\begin{ass}\label{ass_gainop_wd}
The family $\{\gamma_{ij} : i,j \in \N\}$ is pointwise equicontinuous. That is, for every $r \geq 0$ and $\eps>0$, there exists $\delta = \delta(r,\eps) > 0$ such that $|r - \tilde{r}| \leq \delta$, $\tilde{r} \in \R_+$, implies $|\gamma_{ij}(r) - \gamma_{ij}(\tilde{r})| \leq \eps$ for all $i,j\in\N$.
\end{ass}

Additionally, we make the following assumption on the external gains.%

\begin{ass}\label{ass_gammaiu_max}
There is $\gamma_{\max}^u \in \K$ such that $\gamma_{iu} \leq \gamma_{\max}^u$ for all $i \in \N$.%
\end{ass}

We now introduce the concept of a \emph{path of strict decay} (for the gain operator $\Gamma$) which is of crucial importance in the construction of an ISS Lyapunov function for the interconnected system.

\begin{definition}\label{def_omega_path}
A mapping $\sigma:\R_+ \rightarrow \ell_{\infty}^+$ is called a \emph{path of strict decay (for $\Gamma$)}, if the following properties hold:%
\begin{enumerate}
\item[(i)] There exists a function $\rho \in \Kinf$ such that%
\begin{equation*}
  \Gamma(\sigma(r)) \leq (\id + \rho)^{-1} \circ \sigma(r) \mbox{\quad for all\ } r \geq 0,%
\end{equation*}
where $(\id + \rho)^{-1}$ is applied componentwise.%
\item[(ii)] There exist $\sigma_{\min},\sigma_{\max} \in \Kinf$ satisfying%
\begin{equation*}
  \sigma_{\min} \leq \sigma_i \leq \sigma_{\max} \mbox{\quad for all\ } i \in \N.%
\end{equation*}
\item[(iii)] Each component function $\sigma_i = \pi_i \circ \sigma$, $i\in\N$, is a $\Kinf$-function.%
\item[(iv)] For every compact interval $K \subset (0,\infty)$, there exist $0 < c \leq C < \infty$ such that for all $r_1,r_2 \in K$ and $i \in \N$%
\begin{equation*}
  c|r_1 - r_2| \leq |\sigma_i^{-1}(r_1) - \sigma_i^{-1}(r_2)| \leq C|r_1 - r_2|.%
\end{equation*}
\end{enumerate}
\end{definition}

Recall that $V:X\to\R_+$ is called locally Lipschitz continuous on $S \subset X$, i.e. for each $x_0\in S$ there is $\delta>0$ and $L>0$ such that for all $x_1,x_2 \in \{x\in X:\|x-x_0\|\leq \delta\}$ it holds that
\[
|V(x_1)-V(x_2)|\leq L\|x_1-x_2\|_X.
\]

Now we can formulate a Lyapunov-based small-gain result for ISS of infinite interconnections due to \cite{KMZ23}:
\begin{theorem}
\label{thm_mainres}
Consider the (well-posed) interconnected system $\Sigma$, composed of subsystems $\Sigma_i$, $i\in\N$, and let each of the subsystems possess an ISS Lyapunov function as in Assumption \ref{ass_subsystem_iss}.
Let further
\begin{enumerate}
\item There exists a path $\sigma:\R_+ \rightarrow \ell_{\infty}^+$ of strict decay for the gain operator $\Gamma$, defined via the internal gains $\gamma_{ij}$.%
\item For each $R>0$, there is a constant $L(R)>0$ such that%
\begin{equation}\label{eq_lyap_lipschitz}
  |V_i(x_i) - V_i(y_i)| \leq L(R)|x_i - y_i|%
\end{equation}
for all $i\in\N$ and $x_i,y_i \in B_R(0) \subset \R^{n_i}$.%
\end{enumerate}
Then $\Sigma$ is ISS and an ISS Lyapunov function for $\Sigma$ is given by%
\begin{equation}\label{eq_def_V}
  V(x) := \sup_{i\in\N} \sigma_i^{-1}(V_i(x_i)) \mbox{\quad for all\ } x \in X.%
\end{equation}
Moreover, $V$ is locally Lipschitz continuous on $X \setminus \{0\}$. 
\end{theorem}


\begin{ack}[Acknowledgments]

A. Mironchenko is supported by the German Research Foundation (DFG) by Heisenberg project MI~1886/3-1.
\end{ack}

\seealso{article title article title. You can consider mentioning another contribution by A. Mironchenko to this Encyclopedia.}

\bibliographystyle{Harvard}
\bibliography{Mir_LitList_NoMir,MyPublications}

\end{document}